\documentclass[leqno,twoside,11pt]{amsart}

\usepackage{latexsym,esint,url}
\usepackage{psfig}

\theoremstyle{plain}

\newtheorem{theorem}{Theorem}[section]

\newtheorem{conjecture}[theorem]{Conjecture}
\theoremstyle{definition}

\numberwithin{equation}{section}
\numberwithin{figure}{section}

\newcommand{\e}{\epsilon}
\newcommand{\R}{\mathbb {R}}
\newcommand{\Div}{\mbox{div}}
\newcommand{\ds}{\displaystyle}

\begin{document}

\title[infinite hierarchy of shell models]
{The infinite hierarchy of elastic shell models: \\
some recent results and a conjecture}
\author{Marta Lewicka and Reza Pakzad}
\address{Marta Lewicka, University of Minnesota, Department of Mathematics, 
206 Church St. S.E., Minneapolis, MN 55455}
\address{Reza Pakzad, University of Pittsburgh, Department of Mathematics, 
139 University Place, Pittsburgh, PA 15260}
\email{lewicka@math.umn.edu, pakzad@pitt.edu}

\date{\today}
\begin{abstract} 
We summarize some recent results of the authors and their collaborators, 
regarding the derivation of thin elastic shell models (for shells with mid-surface
of arbitrary geometry) from the variational theory of $3$d nonlinear elasticity. 
We also formulate a conjecture on the form and validity of infinitely 
many limiting $2$d models, each corresponding to its proper scaling range of 
the body forces in terms of the shell thickness.
\end{abstract}

\maketitle
\tableofcontents

\section{Introduction} 
 
Elastic materials exhibit qualitatively different responses to different kinematic 
boundary conditions or body forces. A sheet of paper may crumple under 
compressive forces, but it shows a more rigid behavior in a milder regime.
A cylinder buckles in presence of axial loads.  A clamped convex shell enjoys 
great resistance to bending and stretching, but if a hole is pierced into it, 
the whole structure might easily collapse. 
Growing tissues, such as leaves, attain non-flat elastic equilibrium 
configurations with non-zero stress, even in the absence of any external forces. 

Such observations gave rise to many interesting questions in the mathematical 
theory of elasticity. Its main goal is to explain these apparently different 
phenomena based on some common mathematical ground.  
Among others, the variational approach to the nonlinear theory has been 
very effective in rigorously deriving models pertaining to different scaling 
regimes of the body forces \cite{FJMhier}.  The strength of this approach 
lies in its ability to predict the appropriate model together with the response
of the plate without any a priori assumptions other than the general principles 
of $3$d nonlinear elasticity. 
 
The purpose of this paper is to introduce some new results and conjectures 
on the variational derivation of shell theories. 
They can be considered as generalizations
of the results in \cite{FJMhier}, justifying a hierarchy of theories for 
nonlinearly elastic plates. This hierarchy corresponds to the scaling of 
the elastic energy in terms of thickness $h$, in the limit as $h\to 0$.
Some of the derived models were absent from the physics and engineering 
literature before. 

\bigskip

\noindent{\bf Acknowledgments.} 
M.L. was partially supported by the NSF grant DMS-0707275. 
and by the Center for Nonlinear Analysis (CNA) under 
the NSF grants 0405343 and 0635983.
R.P. was partially supported by the University of Pittsburgh 
grant CRDF-9003034.

\section{The set-up and a glance at previously known results}

\noindent
{\bf 2.1. Three dimensional nonlinear elasticity and 
the limiting lower dimensional theories.} 
The equations for the balance of linear momentum for the deformation $u= u(t,x)
\in \mathbb{R}^3$ of the reference configuration $\Omega \subset \mathbb{R}^3$ 
of an elastic body with constant temperature and density read \cite{Ba1}: 
\begin{equation}\label{balancelaw}   
\partial_{tt} u - \Div\,  DW (\nabla u) = f, 
\end{equation} where $DW$ is the Piola-Kirchhoff stress tensor, $f$ is the external
body force, and the elastic energy density $W$ is assumed to satisfy 
the following fundamental properties of frame indifference
(with respect to the group of proper rotations $SO(3)$), normalization 
and non-degeneracy:
\begin{equation}\label{energyW}
\begin{split}
\forall F\in \mathbb{R}^{3\times 3} \quad
\forall R\in SO(3) \qquad
& W(RF) = W(F), \quad W(R) = 0,\\
&W(F)\geq c \cdot\mathrm{dist}^2(F, SO(3)), 
\end{split}
\end{equation}
with a uniform constant $c>0$.

The steady state solutions to (\ref{balancelaw}) satisfy the equilibrium equations:
$- \Div\,  DW (\nabla u) = f$ which, expressed in their weak form, yield 
the formal Euler-Lagrange equations for the critical points of 
the total energy functional: 
\begin{equation}\label{total-energy} 
J(u) = \int_\Omega W(\nabla u) - \int_\Omega f u,   
\end{equation} 
defined for deformation $u: \Omega \longrightarrow \R^3$. 
We will refer to the term 
$E(u) =  \ds \int_\Omega W(\nabla u)$ as the elastic energy of the deformation $u$. 

As a first step towards understanding the dynamical problem (\ref{balancelaw})
it is natural to study the minimizers of (\ref{total-energy}), in an appropriate 
function space. The questions regarding existence and regularity 
of these minimizers are vastly considered in the literature. 
However, due to the loss of convexity of $W$, caused by the frame indifference 
assumption, these problems cannot be dealt with the usual techniques in the 
calculus of variations; see \cite{Ba1} for a review of results and open problems. 

One advantageous direction of research has been to restrict the attention to
domains $\Omega$ which are thin in one or two directions, and hence practically 
reduce the theory to a $2$d or $1$d problem. Indeed, the derivation of 
lower dimensional models for thin structures (such as membranes, shells, or beams) 
has been one of the fundamental questions since the beginning of research 
in elasticity \cite{Love}. The classical approach is to propose a formal asymptotic 
expansion for the solutions  (in other words an {\it Ansatz}) 
and derive the corresponding limiting theory by considering the first terms of the 
$3$d equations under this expansion \cite{ciarbookvol}. The more rigorous 
variational approach of {\em $\Gamma$-convergence} was more recently applied by 
LeDret and Raoult \cite{LeD-Rao} in this context, and then significantly furthered 
by Friesecke, James and M\"uller \cite{FJMhier}, leading to the derivation of 
a hierarchy of limiting plate theories. Among other features, it provided 
a rigorous justification of convergence of minimizers of (\ref{total-energy}) 
to minimizers of suitable lower dimensional limit energies. 

\medskip

\noindent
{\bf 2.2. $\Gamma$-convergence.}
Recall \cite{dalmaso} that a sequence of functionals 
$F_n: X \longrightarrow [-\infty, +\infty]$ defined on a metric space $X$, 
$\Gamma$-converges to the limit functional 
$F: X\longrightarrow [-\infty, +\infty]$ whenever:
\begin{itemize}
\item[{(i)}] (the $\Gamma$-liminf inequality) For any sequence $x_n \to x$ in $X$, 
one has $F(x) \le \liminf_{n \to \infty} F_n(x_n)$.
\item[{(ii)}] (the $\Gamma$-limsup inequality) For any $x\in X$, there exits 
a sequence $x_n$ (called a {\it recovery sequence}) converging to $x\in X$, 
such that $\limsup_{n\to \infty} F_n(x_n) \le F(x)$. 
\end{itemize}  
It is straightforward that $\lim_{n\to \infty} F_n(x_n) = F(x)$ for any recovery 
sequence $x_n \to x$.   

When $X$ is only a topological space, the definition 
of $\Gamma$-convergence involves, naturally, 
systems of neighborhoods rather than sequences. However, when the functionals $F_n$ 
are equi-coercive and $X$ is a reflexive Banach space equipped with weak topology, 
one can still use (i) and (ii) above (for weakly converging sequences), 
as an equivalent version of this definition.

A fundamental consequence of the above definition is the following. If $x_n$ 
is a sequence of approximate minimizers of $F_n$ in $X$:
$$ \lim_{n\to \infty} \left\{F_n(x_n) - \inf_X F_n\right\} =0,$$ 
and if $x_n \to x$, then $x$ is a minimizer of $F$. 
In turn, any recovery sequence associated to a minimizer of $F$ is an approximate 
minimizing sequence for $F_n$.  The convergence of (a subsequence of) $x_n$
is usually independently established through a compactness argument.
 
\medskip

\noindent
{\bf 2.3. A glance at previously known results.}
Let $S$ be a $2$d surface embedded in $\mathbb{R}^3$, which is compact, connected, 
oriented, of class $\mathcal{C}^{1,1}$ and whose boundary $\partial S$ is the union
of finitely many (possibly none) Lipschitz curves.
By $\vec n$ we denote the unit normal vector to $S$, and  
$\pi :S^{h_0} \to S$ is the usual orthogonal projection of the tubular neighborhood
onto $S$.

Consider a family $\{S^h\}_{h>0}$ of thin shells of thickness $h$ around $S$:
$$S^h = \{z=x + t\vec n(x); ~ x\in S, ~ -h/2< t < h/2\}, \qquad 0<h<h_0, $$
The elastic energy per unit thickness of a deformation 
$u\in W^{1,2}(S^h, \mathbb{R}^3)$ is given by:
\begin{equation}\label{elastic-En}
E^h(u) = \frac{1}{h}\int_{S^h} W(\nabla u),
\end{equation}
On above the properties 
where the stored-energy density $W:\mathbb{R}^{3\times 3}\longrightarrow [0,\infty]$
is assumed to satisfy (\ref{energyW}) and to be $\mathcal{C}^2$ regular 
in some open neighborhood of $SO(3)$.

In presence of applied forces $f^h\in L^2(S^h, \mathbb{R}^3)$, the (scaled) total
energy reads:
\begin{equation}\label{total-intro}
J^h(u) = E^h(u) - \frac{1}{h}\int_{S^h} f^hu.
\end{equation} 
It can be shown that if the forces $f^h$ scale like $h^\alpha$, then the 
elastic energy $E^h(u^h)$ at (approximate) minimizers $u^h$ of $J^h$
scale like $h^\beta$, where $\beta= \alpha$ if $0 \le \alpha \le 2$ and $\beta
= 2\alpha -2$ if $\alpha > 2$. The main part of the analysis consists 
therefore of identifying the $\Gamma$-limit $I_\beta$ of the energies
$h^{-\beta} E^h(\cdot, S^h)$  as $h\to 0$, for a given scaling $\beta\geq 0$.  
No a priori assumptions are made on the form of the deformations $u^h$
in this context.

In the case when $S$ is a subset of $\mathbb R^2$ (i.e. a plate), 
such $\Gamma$-convergence was first established for $\beta = 0$ \cite{LR1},
and  later \cite{FJMgeo, FJMhier} for all $\beta \geq 2$. 
This last scaling regime corresponds to a rigid  behavior of the elastic material, 
since the limiting admissible deformations are isometric immersions 
(if $ \beta =2$) or infinitesimal isometries (if $\beta>2$) of the mid-plate $S$.
One particular case is  $\beta=4$, where the derived limiting theory turns 
out to be the  von K\'arm\'an theory \cite{karman}. 
A totally clamped plate exhibits a very rigid behavior already for $\beta>0$ 
\cite{CMM}.
In case $0 < \beta < 5/3$, the $\Gamma$-convergence was recently 
obtained in \cite{CM05}, while the regime $5/3 \le \beta<2$ remains open and 
is conjectured to be relevant for crumpling of elastic sheets \cite{shankar}. 
 
Much less is known in the general case when $S$ is a surface of arbitrary geometry.
The first result in \cite{LeD-Rao} relates to scaling $\beta=0$ and models 
{\em membrane shells}: the limit $I_0$ depends only on the stretching 
and shearing produced by the deformation on $S$. 
Another study \cite{FJMM_cr} analyzed the case $\beta=2$, corresponding to 
a {\em flexural shell model} \cite{ciarbookvol}, or a geometrically nonlinear 
purely bending theory, where the only admissible deformations are 
isometric immersions, that is those preserving the metric on $S$ (see section 2). 
The energy $I_2$ depends then on the change of curvature produced 
by the deformation. 
 
All the above mentioned theories should be put in contrast with a large body 
of literature, devoted to derivations starting from $3$d {\em linear} elasticity 
(see \cite{ciarbookvol} and references therein). In the present setting one 
allows for large deformations, i.e. not necessarily close to a rigid motion.
The basic assumption of the linear elasticity is not taken for granted 
in our context.

\section{The Kirchhoff theory for shells: $\beta=2$ and arbitrary $S$}

The limiting theory for $\beta=2$ is precisely described in the following result:
\begin{theorem}\cite{FJMM_cr}\label{kirchhoff} 
(a) {\bf Compactness and the $\Gamma$-liminf inequality.}  
Let $u^h \in W^{1,2}(S^h, \R^3)$ be a sequence of deformations such that 
$E^h(u^h)/h^2$ is uniformly bounded. 
Then there exists a sequence $c^h \in \R^3$ such that the rescaled deformations:
\begin{equation*}\label{rescale}
y^h (x + t \vec n) = u^h (x + t {h}/{h_0}\vec n) - c^h: S^{h_0}\longrightarrow 
\mathbb{R}^3,
\end{equation*} 
converge (up to a subsequence) in $W^{1,2}$ to $y \circ \pi$,
where $y\in W^{2,2} (S, \R^3)$  and it satisfies:
\begin{equation}\label{isom}
(\nabla y)^T \nabla y = \mathrm{Id} \qquad \mbox{a.e. in } S.
\end{equation}
Moreover: 
$$I_2(y)=  \int_S \mathcal Q_2(x, \Pi(y) - \Pi) \le 
\liminf_{h\to 0} \frac 1{h^2}E^h (u^h).$$ 
(b) {\bf The recovery sequence and the $\Gamma$-limsup inequality.} 
Given any isometric immersion $y\in W^{2,2}(S, \R^3)$ satisfying (\ref{isom}), 
there exists a sequence $u^h \in W^{1,2}(S^h,\R^3)$
such that the rescaled deformations $y^h(x+ t\vec n) = u^h(x+th/h_0\vec n)$
converge to $y \circ \pi$  in $W^{1,2}$ and:
\begin{equation*}\label{limsupkir}
I_2(y) \ge \limsup_{h\to 0} \frac 1{h^2}E^h (u^h).
\end{equation*} 
\end{theorem}
In the definition of the limit functional $I_2$,
the quadratic forms $\mathcal{Q}_2(x,\cdot)$ are defined as follows:
\begin{equation}\label{Q2}
\mathcal{Q}_2(x, F_{tan}) = \min\{\mathcal{Q}_3(\tilde F); 
~~ (\tilde F - F)_{tan} = 0\}, 
\qquad \mathcal{Q}_3(F) = D^2 W(\mbox{Id})(F,F).
\end{equation}
The form $\mathcal{Q}_3$ is defined for all $F\in\mathbb{R}^{3\times 3}$, 
while $\mathcal{Q}_2(x,\cdot)$,  for a given $x\in S$ is defined on tangential 
minors $F_{tan}$ of such matrices.   Recall that the tangent 
space to $SO(3)$ at $\mbox{Id}$ is $so(3)$. As a consequence,  
both forms depend only on the symmetric parts of their arguments
and are positive definite on the space of symmetric matrices \cite{FJMgeo}.

The functional $I_2(y)$ measures the total change of curvature (bending) 
induced by the deformation $y$ of the mid-surface $S$. In the form of 
the integrand $\Pi$ denotes the shape operator on $S$, 
while $\Pi(y)$ is the pull back of the 
shape operator of the surface $y(S)$ under $y$.
For any orthonormal tangent frame $\tau, \eta \in T_x S$ there holds:
$$\eta \cdot \Pi \tau = \eta \cdot  \partial_\tau\vec n   \quad \mbox{and} \quad
\eta \cdot \Pi(y) \tau = \eta \cdot \partial_\tau \vec N,$$ 
where $\vec N: S\to \R^3$ is the unit normal to $y(S)$:
$\vec N (x) = \partial_\tau y \times \partial_\eta y.$   

\section{The von-K\'arm\'an theory for shells: $\beta=4$ and arbitrary $S$}

For the range of scalings $\beta>2$ a rigidity argument 
\cite{FJMgeo, FJMhier, lemopa1} shows that the admissible
deformations $u$ are only those which are close to a rigid motion
$R$ and whose first order term in the expansion of $u-R$ with respect to $h$
is given by $RV$. The displacement field $V$ is an element of the 
class $\mathcal V_1$ of {\em infinitesimal isometries} on $S$ \cite{Spivak}. 
The space $\mathcal{V}_1$ consists of vector fields
$V\in W^{2,2}(S,\mathbb{R}^3)$ for whom there 
exists  a matrix field $A\in W^{1,2}(S,\mathbb{R}^{3\times 3})$ so that:
\begin{equation}\label{Adef-intro} 
\partial_\tau V(x) = A(x)\tau \quad \mbox{and} \quad  A(x)^T= -A(x) \qquad 
{\rm{a.e.}} \,\, x\in S \quad \forall \tau\in T_x S.
\end{equation}       
In other terms, $V$ is a (first order) infinitesimal isometry if the change 
of metric induced by the deformation $\mbox{id} + \e V$ is at most of order $\e^2$.
\begin{figure}[h] 
\centerline{\psfig{figure=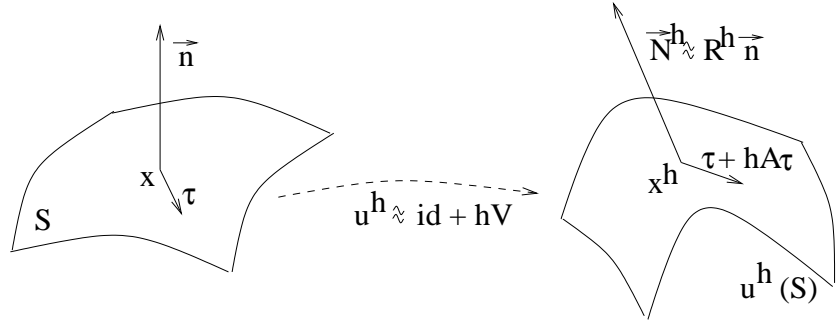,width=11cm,angle=0}}
\caption{The mid-surface $S$ and its deformation.}
\label{figura1}
\end{figure}

For $\beta=4$ the $\Gamma$-limit turns out to be the generalization of 
the von K\'arm\'an functional \cite{FJMhier} to shells, 
and it consists of two terms: 
\begin{equation*}\label{vonKarman}
I_4(V,B_{tan})= \frac{1}{2} 
\int_S \mathcal{Q}_2\left(x,B_{tan} - \frac{1}{2} (A^2)_{tan}\right)
+ \frac{1}{24} \int_S \mathcal{Q}_2\left(x,(\nabla(A\vec n) - A\Pi)_{tan}\right).
\end{equation*}
The quadratic form $\mathcal {Q}_2$ is defined as in (\ref{Q2}) and $A$ is as in
(\ref{Adef-intro}).  The second term above measures bending, that is 
the first order change in the second fundamental 
form of $S$, produced by $V$. The first term measures stretching, that is the 
second order change in the metric of $S$.
It involves a symmetric matrix field $B_{tan}$ belonging to the 
{\em finite strain space}: 
$$ \mathcal{B} = \mbox{cl}_{L^2}\Big\{\mathrm{sym }\nabla w; 
~ w\in W^{1,2}(S,\mathbb{R}^3)\Big\}. $$ 
The space $\mathcal{B}$ emerges as well in the context of  linear elasticity
and ill-inhibited surfaces \cite{sanchez, GSP}.
\begin{theorem}\cite{lemopa1}\label{thmvonkarman}
(a) Let $u^h\in W^{1,2}(S^h,\mathbb{R}^3)$ be a sequence of deformations whose 
scaled energies $E^h(u^h)/h^4$ are uniformly bounded.
Then there exist a sequence $Q^h\in SO(3)$ and $c^h\in\mathbb{R}^3$ 
such that for the normalized rescaled deformations:
\begin{equation}\label{norm2}
y^h(x+t\vec{n}) = Q^h u^h(x+h/h_0 t\vec{n}) - c^h:S^{h_0}\longrightarrow \R^3
\end{equation}
the following holds.
\begin{enumerate}
\item[(i)] $y^h$ converge in $W^{1,2}(S^{h_0})$ to $\pi$.
\item[(ii)]  The scaled average displacements:
\begin{equation*}\label{Vh-intro}
V^h(x) = \frac{1}{h} \fint_{-h_0/2}^{h_0/2}
y^h(x+t\vec{n}) - x ~\mathrm{d}t
\end{equation*} 
converge (up to a subsequence) in $W^{1,2}(S)$ to some $V\in \mathcal{V}_1$.
\item[(iii)] The scaled strains $\frac{1}{h} ~\mathrm{ sym}\nabla V^h$ 
converge weakly
in $L^2$ to a symmetric matrix field $B_{tan} \in \mathcal {B}$. 
\item[(iv)] $ I_4(V, B_{tan}) \le \liminf_{h\to 0} {1}/{h^4} E^h(u^h).$
\end{enumerate}
(b) For every $V\in\mathcal{V}_1$ and $B_{tan} \in \mathcal {B}$ there exists a sequence
of deformations $u^h\in W^{1,2}(S^h,{\mathbb R}^3)$ such that:
\begin{itemize}
\item[(i)] the rescaled deformations $y^h(x+t\vec n)=u^h(x+th/h_0 \vec n)$ 
converge in $W^{1,2}(S^{h_0})$ to $\pi$.
\item[(ii)] the scaled average displacements $V^h$ given above converge 
in $W^{1,2}(S)$ to $V$.
\item[(iii)] the scaled linearized strains $\frac{1}{h} ~\mathrm{ sym}\nabla V^h$ 
converge weakly in $L^2$ to  $B_{tan}$. 
\item[(iv)] $  I_4(V, B_{tan}) \ge \limsup_{h\to 0} {1}/{h^4} E^h(u^h) $.
\end{itemize}
\end{theorem}
The special case of this theorem for plates, that is when $S \subset \R^2$, 
was already proved in \cite{FJMhier}. It can be shown that for a flat surface, 
the infinitesimal isometries coincide essentially with the out-of-plane 
displacements. Also, the space $\mathcal B$ becomes then the set of all linearized 
strains associated with the in-plane displacements in $S$, so the functional 
$I_4$ can be written directly in terms of an out-of-plane and an in-plane 
displacement. The Euler-Lagrange equations derived from this limit functional 
lead to the von-K\'arm\'an equations \cite{karman}. 

\section{The linear theory for shells:  $\beta>4$ and arbitrary $S$;
$\beta=4$ and approximately robust $S$ }

It was shown in \cite{lemopa1} that for a certain class of surfaces, 
referred to as {\it approximately robust surfaces}, the limiting 
theory for $\beta=4$ reduces to the purely linear bending functional:
\begin{equation}\label{linearbend}
I_{lin}(V)= \frac{1}{24} \int_S 
\mathcal{Q}_2\left(x,(\nabla(A\vec{n}) - A\Pi)_{tan}\right)~\mbox{d}x
\qquad \forall V\in\mathcal{V}_1,
\end{equation} 
This class of surfaces is given by the property that any 
first order infinitesimal isometry $V\in \mathcal {V}_1$ can be modified to 
be arbitrarily close to a $W^{1,2}$ second order isometry:
$$\mathcal {\tilde V}_2= \{ V\in \mathcal {V}_1;~ (A^2)_{tan} \in {\mathcal B} \}
= \mathcal {V}_1.$$ 
Convex surfaces, surfaces of revolution and developable surfaces  
belong to this class \cite{lemopa1}.

In \cite{lemopa1}, it was also proved that the $\Gamma$-limit of $E^h/h^\beta$ 
for the scaling regime $\beta>4$ is also given by the functional (\ref{linearbend}).
This corresponds to the linear pure bending theory derived in \cite{ciarbookvol} 
from linearized elasticity.  The important qualitative difference between 
this theory and the limiting theory for $\beta = 4$ and an approximately robust 
surface is in the type of convergences one establishes for a sequence $u^h$ 
satisfying $E^h(u^h) \le C h^\beta$. Indeed, if $\beta>4$, the best one can prove 
is the convergence in $W^{1,2}$, up to a subsequence, of the rescaled displacement 
fields:
$$ V^h(x) = \frac{1}{h^{\beta/2-1}} \fint_{-h_0/2}^{h_0/2}
y^h(x+t\vec{n}) - x ~\mathrm{d}t$$ 
to an element $V\in \mathcal {V}_1$. Note the finer rescaling parameter 
$h^{\beta/2-1}$ with respect to (\ref{Vh-intro}).

\section{Intermediate theories for plates and convex shells: $\beta\in (2,4)$ }
 
In paper \cite{lemopa3} we focused on the range of scalings $2<\beta<4$,
looking hence for an intermediate theory between those corresponding 
to $\beta=2$ and $\beta\ge 4$. On one hand, modulo a rigid motion, 
the deformation of the mid-surface must 
look like $\mbox{id} + \e V$, up to its first order of expansion. 
On the other hand, the closer $\beta$ is to $2$, 
the closer the mid-surface deformation must be to an exact isometry of $S$. 
To overcome this apparent disparity between first order infinitesimal
isometries and exact isometries in this context, one is immediately drawn 
to consider higher order infinitesimal isometries which lay somewhat between 
these two categories. This will be the subject of discussion in section
\ref{conjecture}. Another angle of approach, which turns out to be 
useful in special cases, is to study conditions under which, 
given $V \in \mathcal {V}_1$, one can construct an exact isometry of the form 
$\mbox{id} + \e V + \e^2 w_\e$, with equibounded corrections $w_\e$. This is what 
we refer to as a {\em matching property}.

If $S\subset \mathbb R^2$ represents a plate, the above issues have been
answered  in \cite{FJMhier}. In this case:
\begin{itemize}
\item[(i)] The limit displacement $V$ must necessarily belong to the space of
second order infinitesimal isometries:
$ \mathcal {\tilde V}_2= \{ V\in \mathcal {V}_1;~ (A^2)_{tan} \in {\mathcal B} \},$ 
where the matrix field $A$ is as in (\ref{Adef-intro}).  
\item[(ii)] Any Lipschitz second order isometry $V\in\mathcal {\tilde V}_2$
satisfies the matching property. 
\end{itemize}  
Combining these two facts with the density of Lipschitz second order 
infinitesimal isometries in $\mathcal {\tilde V}_2$ for a plate \cite{MuPa2},  
one concludes through the $\Gamma$-convergence arguments that the 
limiting plate theory is given by the functional (\ref{linearbend}) over
$\mathcal {\tilde V}_2$.  Note that, for a plate, $V\in \mathcal {\tilde V}_2$ 
means that there exists an in-plane displacement $w\in W^{1,2} (S, \mathbb R^2)$ 
such that the change of metric due to $\mbox{id} + \e V + \e^2 w$ is of order 
$\e^3$. Also, in this case, an equivalent analytic characterization for 
$V=(V^1, V^2, V^3) \in \mathcal {\tilde V}_2$ is given by 
$(V^1, V^2) =  (-\omega y, \omega x) + (b_1, b_2)$ and: $\det \nabla^2 V^3 =0$. 

\medskip

Towards analyzing more general surfaces $S$, we derived a matching
property and the corresponding density of isometries, for elliptic surfaces.
We say that $S$ is elliptic if its shape operator $\Pi$ is strictly positive 
(or strictly negative) definite up to the boundary:
\begin{equation}\label{elliptic-intro} 
\forall x\in \bar S \quad \forall \tau\in T_xS \qquad 
\frac{1}{C}|\tau|^2 \leq \Big(\Pi(x)\tau\Big)\cdot\tau \leq C|\tau|^2.
\end{equation}  
The novelty here is the fact that for an elliptic surface, 
all sufficiently smooth infinitesimal isometries satisfy the matching property:
\begin{theorem}\label{th_exact-intro} \cite{lemopa3}
Let $S$ be elliptic as in \eqref{elliptic-intro}, homeomorphic to a disk 
and let for some $\alpha >0$,  $S$ and $\partial S$ be of class 
$\mathcal{C}^{3,\alpha}$. Given $V\in\mathcal{V}_1\cap\mathcal{C}^{2,\alpha}(\bar S)$,
there exists a sequence $w_\e:\bar S\longrightarrow \mathbb{R}^3$, equibounded in 
$\mathcal{C}^{2,\alpha}(\bar S)$, and such that for all small $\e>0$ the map
$u_\e = \mathrm{id} +\e V + h^2w_\e$ is an (exact) isometry. 
\end{theorem}  
We apply this result to construct the recovery sequence in the $\Gamma$-limsup 
inequality.  Clearly, Theorem \ref{th_exact-intro} is not sufficient for this 
purpose as the elements of $\mathcal {V}_1$ are only $W^{2,2}$ regular. In 
most $\Gamma$-convergence results, a key step is to prove density of suitably 
regular mappings in the space of mappings admissible for the limit problem. 
Results in this direction, for Sobolev spaces of isometries and infinitesimal
isometries, have been shown and applied in the context of derivation 
of plate theories. The interested reader can refer to \cite{Pak, MuPa2} 
for statements of these density theorems and their applications in 
\cite{FJMhier, CD}.

In general, even though $\mathcal{V}_1$ is a linear space, and assuming $S$ to be  
$\mathcal{C}^\infty$, the usual mollification techniques do not guarantee 
that elements of $\mathcal{V}_1$ can be approximated by smooth 
infinitesimal isometries. An interesting example, discovered by Cohn-Vossen
\cite{Spivak}, is a closed smooth surface of non-negative curvature for which 
$\mathcal{C}^\infty \cap \mathcal{V}_1$ consists only of trivial fields 
$V: S\longrightarrow {\mathbb R}^3$ with constant gradient, 
whereas $\mathcal{C}^2 \cap \mathcal{V}_1$ contains non-trivial infinitesimal 
isometries. Therefore $\mathcal{C}^\infty \cap \mathcal{V}_1$ is not dense 
in $\mathcal{V}_1$ for this surface. We however have:
\begin{theorem}\label{th_density-intro} \cite{lemopa3}
Assume that $S$ is elliptic, homeomorphic to a disk, of class 
$\mathcal{C}^{m+2,\alpha}$ 
up to the boundary and that $\partial S$ is $\mathcal{C}^{m+1,\alpha}$, 
for some $\alpha\in (0,1)$ and an integer $m>0$. 
Then, for every $V\in\mathcal{V}_1$ there exists a sequence 
$V_n\in\mathcal{V}_1\cap\mathcal{C}^{m,\alpha}(\bar S,\mathbb{R}^3) $ 
such that: $$\lim_{n\to\infty} \|V_n - V\|_{W^{2,2}(S)} = 0.$$
\end{theorem} 
Ultimately, and as a consequence of Theorems \ref{th_exact-intro} and 
\ref{th_density-intro}, the main result of \cite{lemopa3} states that for 
elliptic surfaces of sufficient regularity, the $\Gamma$-limit of 
the nonlinear elastic energy (\ref{elastic-En}) for the scaling regime
$2<\beta<4$ (and hence for all $\beta>2$) is still given by the functional 
(\ref{linearbend}) over the linear space ${\mathcal V_1}$:
\begin{theorem}\label{thliminf-intro} \cite{lemopa3}
Let $S$ be as in Theorem \ref{th_exact-intro} and let $2<\beta <4$.

\noindent (a) Assume that for a sequence of deformations 
$u^h\in W^{1,2}(S^h,\mathbb{R}^3)$ their 
scaled energies $E^h(u^h)/ h^\beta$ are uniformly bounded.
Then there exist a sequence $Q^h\in SO(3)$ and $c^h\in\mathbb{R}^3$ 
such that for the normalized rescaled deformations in (\ref{norm2}) 
the following holds.
\begin{enumerate}
\item[(i)] $y^h$ converge in $W^{1,2}(S^{h_0})$ to $\pi$.
\item[(ii)]  The scaled average displacements:
\begin{equation*}\label{Vh-indtro}
V^h(x) = \frac{1}{h^{\beta/2 -1}} \fint_{-h_0/2}^{h_0/2}
y^h(x+t\vec{n}) - x ~\mathrm{d}t
\end{equation*} 
converge (up to a subsequence) in $W^{1,2}(S)$ to some $V\in \mathcal{V}_1$.
\item[(iii)] $I_{lin}(V) \le \liminf_{h\to 0} {1}/{h^\beta} E^h(u^h).$
\end{enumerate}
(b) For every $V\in\mathcal{V}_1$ there exists a sequence
of deformations $u^h\in W^{1,2}(S^h,{\mathbb R}^3)$ such that:
\begin{itemize}
\item[(i)] the rescaled deformations $y^h(x+t\vec n)=u^h(x+th/h_0 \vec n)$ 
converge in $W^{1,2}(S^{h_0})$ to $\pi$.
\item[(ii)] the scaled average displacements $V^h$ given above converge 
in $W^{1,2}(S)$ to $V$.
\item[(iii)] $  I_{lin}(V) \ge \limsup_{h\to 0} {1}/{h^\beta} E^h(u^h) $.
\end{itemize}
\end{theorem} 
One can actually prove that for $2<\beta<4$ and surface $S$ of arbitrary geometry,
the part (a) or Theorem \ref{thliminf-intro} remains valid, and moreover
$1/h^{\beta/2-1} ~\mathrm{sym }~\nabla V^h$ converge
(up to a subsequence) in $L^{2}(S)$ to ${1}/{2}(A^2)_{tan}$, where $A$ is
related to $V$ by (\ref{Adef-intro}).
The novelty with respect to the equivalent result for $\beta =4$
in the first part of Theorem \ref{thmvonkarman} is the constraint 
$V\in \mathcal {\tilde V}_2$. If $S$ is an elliptic surface 
of sufficient regularity, the set $\mathcal B$ coincides with the whole space 
$L^2_{sym}(S, {\mathbb R}^{2\times2})$ \cite{lemopa1}, hence the constraint 
is automatically satisfied for all $V\in \mathcal {V}_1$.  
In the general case where $S$ is an arbitrary surface, a characterization of 
this constraint and the exact form of $\mathcal B$ may be complicated.

\medskip

The case of the scaling range $2<\beta<4$ is still open for general shells.  
The following section is dedicated to the presentation of a conjecture on 
this problem, stating that other constraints, similar to the inclusion
$V\in \mathcal {\tilde V}_2$, should be present for values of $\beta$ closer to $2$.
Heuristically, the closer $\beta$ is to $2$, we expect $V$ to be an
infinitesimal isometry of higher order.

\section{A conjecture on the infinite hierarchy of shell models}\label{conjecture}

If the deformations $u^h$ of $S^h$ satisfy a simplified version of Kirchhoff-Love
assumption:
$$u^h(x+t\vec n) = u^h(x) + t \vec N^h(x),$$
vector $\vec N^h$ being the unit normal to the image surface $u^h(S)$,
then formal calculations show that :
\begin{equation}\label{enexpan}
E^h(u^h)\approx \int_S |\delta g_S|^2 + h^2\int_S |\delta \Pi_S|^2.
\end{equation}
Here by $\delta g_S$ and $\delta \Pi_S$ we denote, respectively, the change 
in the metric (first fundamental form) and in the shape operator (second fundamental 
form), between the surface $u^h(S)$ and the reference mid-surface $S$.  
For a more rigorous treatment of this observation see e.g. \cite{CM05}.
The two terms in (\ref{enexpan}) correspond to the stretching and bending energies,
and the factor $h^2$ in the bending term points to the fact that a shell undergoes 
bending more easily than stretching. 
For a plate, the latter energy is known in the literature as the F\"oppl-von K\'arm\'an
functional.

Another useful observation is that for minimizers $u^h$, the energy should be 
distributed equi-partedly between the stretching and bending terms.  
When $E^h(u^h) \approx h^2$, then equating the order of both terms in (\ref{enexpan}) 
we obtain in the limit of $h\to 0$: 
$$\int_S |\delta g_S|^2 \approx 0 \quad \mbox{ and } \quad 
\frac{1}{h^2}E^h(u^h)\approx \int_S |\delta \Pi_S|^2.$$ 
This indeed corresponds to the Kirchhoff model, as in Theorem \ref{kirchhoff} 
\cite{FJMM_cr}, where the limiting energy $I_2$ is given by the bending term 
(measuring the change in the second fundamental form) under the constraint 
of zero stretching: $\delta g_S=0$. As we have seen, the limiting deformation $u$
must be an isometry $(\nabla u)^T\nabla u = \mbox{Id}$ and hence preserve the metric.

\medskip

To discuss higher energy scalings, assume that:
\begin{equation}\label{enscal}
E^h(u^h)\approx h^\beta, \qquad \beta>2.
\end{equation}
Then, as mentioned before, by the rigidity estimate \cite{FJMgeo}, the restrictions 
of $u^h$ to $S$ have, modulo appropriate rigid motions, the following expansions:
$$u^h_{|S} = \mbox{id} + \sum_{i=1}^\infty \epsilon^i w_i.$$
Thus, $\delta\Pi$ is of the order $\epsilon$ and after equating the order of
the bending term in (\ref{enexpan}) by (\ref{enscal}), we arrive at:
$h^2 \epsilon^2 = h^\beta$, that is:
\begin{equation}\label{epbeta}
\epsilon = h^{\beta/2 - 1}.
\end{equation}
On the other hand, the stretching term has the form:
$\delta g_S = (\nabla u^h)\nabla u^h - \mbox{Id} 
= \sum_{i=1}^\infty \epsilon^i A_i$, with:
$$ A_i = \sum_{j+k =i} \mbox{sym} \Big( (\nabla w_j)^T \nabla w_k \Big ), $$ 
indicating the  $i$-th order change of metric. 
Taking into account (\ref{enscal}), this yields: 
$\epsilon^{2i} \int_S |A_i|^2  \approx h^{\beta}$, and so in view of (\ref{epbeta}): 
$\|A_i\|_{L^2}^2 \approx h^{\beta- i(\beta-2)}$. 
A first consequence is that $A_1$ must vanish in the limit as $h\to 0$, that is 
the limiting deformation is a first order infinitesimal isometry. For $i>1$, we observe 
that $\|A_i\|_{L^2}^2\approx h^{(i-1)(\beta_i - \beta)}$, where:
$$\beta_i = 2 + \frac{2}{i-1}.$$ 
We conclude that if $\beta < \beta_N$, then $\|A_i\|_{L^2} \approx 0$ for $i\leq N$, 
and if $\beta = \beta_N$, then $\|A_N\|_{L^2} = O(1)$.  
The study of the asymptotic behavior of the energy $1/h^\beta E^h$ 
leads us hence to the following conjecture. 
\begin{conjecture}
The limiting theory of thin shells with midsurface $S$, under the elastic
energy scaling $\beta>2$ as in (\ref{enscal}) is given by the following
functional $I_\beta$ below, defined on the space $\mathcal{V}_N$
of $N$-th order infinitesimal isometries, where:
$$ \beta\in [\beta_{N+1}, \beta_N).$$
The space $\mathcal{V}_N$ is identified with the space 
of $N$-tuples $(V_1,\ldots, V_N)$ of displacements 
$V_i:S\longrightarrow \mathbb{R}^3$ (having appropriate regularity),
such that the deformations of $S$:
$$u_\epsilon = \mathrm{id} + \sum_{i=1}^N \epsilon^i V_i$$
preserve its metric up to order $\epsilon^N$.
We have: 
\begin{itemize}
\item[(i)] When $\beta=\beta_{N+1}$ then 
$I_\beta = \int_S \mathcal {Q}_2 \left (x, \delta_{N+1}g_S\right)
+ \int_S \mathcal {Q}_2 \left (x, \delta_1 \Pi_S \right)$, 
where $\delta_{N+1}g_S$ is the change of metric
on $S$ of the order $\epsilon^{N+1}$, generated by the family of deformations 
$u_\epsilon$ and $\delta_1 \Pi_S$ is the corresponding first order change in the second 
fundamental form.
\item[(ii)] When $\beta\in (\beta_{N+1}, \beta_N)$ then
$I_\beta = \int_S \mathcal {Q}_2 \left (x, \delta_1 \Pi_S \right )$.
\item[(iii)] The constraint of $N$-th order infinitesimal isometry $\mathcal{V}_N$ 
may be relaxed to that of $\mathcal{V}_M$, $M<N$, if $S$ has the following 
matching property. For every $(V_1,\ldots V_M)\in\mathcal{V}_M$ there exist
sequences of corrections $V_{M+1}^\epsilon,\ldots  V_{N}^\epsilon$, uniformly 
bounded in $\epsilon$, such that:
$$\tilde u_\epsilon = \mathrm{id} + \sum_{i=1}^M \epsilon^i V_i
+ \sum_{i=M+1}^N\epsilon^i V_i^\epsilon$$
preserve the metric on $S$ up to order $\epsilon^N$.
\end{itemize}
\end{conjecture}
This conjecture is consistent with the so far established results in \cite{lemopa1} for 
$N=1$ (i.e. $\beta\ge \beta_2 =  4$) and arbitrary surfaces. 
Note that in the case of approximately robust surfaces, any element of 
$\mathcal{V}_1$ can be matched with an element of $\mathcal {V}_2$, and hence the term 
$\int_S |\delta_{2}g_S|^2$ in the limit energy can be dropped. 
The second order infinitesimal isometry constraint $\mathcal{V}_2$ is established 
for all surfaces when $2<\beta <4$. 
In the particular case of plates, any second order isometry in a dense subset of 
$\mathcal {V}_2$, can be matched with an exact isometry \cite{FJMhier}. 
As a consequence, the theory reduces to minimizing the bending energy under 
the second order infinitesimal isometry constraint. 
A similar matching property for elliptic surfaces, this time for elements of 
$\mathcal {V}_1$, is given in \cite{lemopa3} (see Theorems 
\ref{th_exact-intro} and \ref {th_density-intro}). As a consequence, 
for elliptic surfaces, the limiting theory for the whole range $\beta>2$ 
reduces to the linear bending. 

The case $2<\beta<4$ remains open for all other types of surfaces. 
The main difficulty lies in obtaining the appropriate convergences and 
the limiting nonlinear constraints: 
$$ \sum_{j+k =i} \mbox{sym} \Big ( (\nabla V_j)^T \nabla V_k \Big )=0, \qquad 
1 \le i \le N,$$  
for the elements of $\mathcal{V}_N$ when $\beta < \beta_N$. 
The above nonlinearity, implies a rapid loss of Sobolev regularity of $V_i$
as $i$ increases. Moreover, applying methods of \cite{lemopa3} to surfaces changing 
type leads in this context to working with mixed-type PDEs.


\medskip\begin{thebibliography}{9999}

\bibitem{Ba1} J.M. Ball, \textit{Some open problems in elasticity}, 
in: Geometry, dynamics and mechanics (Marsden Festschrift), Springer, (2002),
pp. 3--59.

\bibitem{ciarbookvol} P.G. Ciarlet, \textit{Mathematical Elasticity, Vol I-III
Theory of Shells}, North-Holland, Amsterdam, (2000).

\bibitem{CD} S. Conti and G. Dolzmann, \textit{$\Gamma$-convergence for incompressible 
elastic plates}, to appear in Calc.Var. PDE, (2008).

\bibitem{CM05} S. Conti and F. Maggi, \textit{Confining thin sheets and folding paper},
Arch. Ration. Mech. Anal.  187  (2008),  no. 1, 1--48. 

\bibitem{CMM} S. Conti, F. Maggi and S. M\"uller,
\textit{Rigorous derivation of F\"oppl's theory for clamped elastic membranes 
leads to relaxation}, SIAM J. Math. Anal.  38,  no. 2, 657--680, (2006).

\bibitem{dalmaso} G. Dal Maso, \textit{An introduction to $\Gamma$-convergence}, 
Progress in Nonlinear Differential Equations and their Applications, {\bf 8}, 
Birkh\"auser, MA, (1993).

\bibitem{FJMM_cr} G. Friesecke, R. James, M.G. Mora and S. M\"uller, 
\textit{Derivation of nonlinear bending theory for shells from three-dimensional 
nonlinear elasticity by Gamma-convergence},  C. R. Math. Acad. Sci. Paris,  
{\bf 336}  (2003),  no. 8, 697--702.

\bibitem{FJMgeo} G. Friesecke, R. James and S. M\"uller, 
\textit{A theorem on geometric rigidity and the derivation of nonlinear 
plate theory from three dimensional elasticity}, Comm. Pure. Appl. Math., 
{\bf 55} (2002), 1461--1506.  

\bibitem{FJMhier} G. Friesecke, R. James and S. M\"uller, \textit{A hierarchy 
of plate models derived from nonlinear elasticity by gamma-convergence}, 
Arch. Ration. Mech. Anal.,  {\bf 180}  (2006),  no. 2, 183--236.

\bibitem{GSP} G. Geymonat and {\'E}. Sanchez-Palencia, \textit{On the rigidity 
of certain surfaces with folds and applications to shell theory},  
Arch. Ration. Mech. Anal., {\bf 129}  (1995),  no. 1, 11--45. 


\bibitem{karman} T. von K\'arm\'an, \textit{Festigkeitsprobleme im Maschinenbau}, 
in Encyclop\"adie der Mathematischen Wissenschaften. Vol. IV/4, pp. 311-385, Leipzig,
1910.

\bibitem{LR1} H. LeDret and A. Raoult, The nonlinear membrane model as a variational 
limit of nonlinear three-dimensional elasticity. \textit{J. Math. Pures Appl.\/} 
\textbf{73} (1995), 549--578.

\bibitem{LeD-Rao} H. LeDret and A. Raoult, \textit{The membrane shell model 
in nonlinear elasticity: a variational asymptotic derivation}, J. Nonlinear Sci., 
{\bf 6} (1996), 59--84.

\bibitem{le} M. Lewicka, \textit{A note on convergence of low energy critical points 
     of nonlinear elasticity functionals, for thin shells of arbitrary geometry},
     submitted.

\bibitem{lemopa1} M. Lewicka, M.G. Mora and M.R. Pakzad, \textit{Shell theories 
arising as low energy $\Gamma$-limit of 3d nonlinear elasticity}, accepted in 
    Annali della Scuola Normale Superiore di Pisa (2009)

\bibitem{lemopa2} M. Lewicka, M.G. Mora and M.R. Pakzad,
    \textit{A nonlinear theory for shells with slowly varying thickness}, 
    C.R. Acad. Sci. Paris, Ser I, {\bf 347}, (2009), 211--216 .

\bibitem{lemopa3} M. Lewicka, M.G. Mora and M. Pakzad,
    \textit{The matching property of infinitesimal isometries on elliptic surfaces and
    elasticity of thin shells}, submitted.

\bibitem{lepa} M. Lewicka and M. Pakzad,
    \textit{Scaling laws for non-Euclidean plates and the 
    $W^{2,2}$ isometric immersions of Riemannian metrics}, submitted.  

\bibitem{LewMul} M. Lewicka and S. M\"uller, \textit{The uniform Korn-Poincar\'e 
inequality in thin domains}, submitted.

\bibitem{Love} A.E.H. Love, \textit{A treatise on the mathematical theory 
of elasticity}, 4th ed.\ Cambridge University Press, Cambridge (1927).

\bibitem{MoSc} M.G. Mora and L. Scardia, \textit{Convergence of equilibria 
of thin elastic plates under physical growth conditions for
the energy density}, preprint (2008).

\bibitem{MuPa2} S. M\"uller and M. R. Pakzad, \textit{Regularity properties of 
isometric immersions},  Math. Z. {\bf 251},  no. 2, 313--331 (2005).  

\bibitem{MuPa} S. M\"uller and M. R. Pakzad, \textit{Convergence of equilibria 
of thin elastic plates -- the von K\'arm\'an case}, Comm. Partial Differential 
Equations {\bf 33},  no. 4-6, 1018--1032, (2008).

\bibitem{Pak} M. R. Pakzad, \textit{On the Sobolev space of isometric immersions}, 
J. Differential Geom. {\bf 66},  no. 1, 47--69 (2004). 

\bibitem{sanchez} {\'E}. Sanchez-Palencia,
\textit{Statique et dynamique des coques minces. II. Cas de flexion pure inhibe\'e. 
Approximation membranaire.}  C. R. Acad. Sci. Paris S\'er. I Math.  {\bf 309}  (1989),  
no. 7, 531--537. 
 
\bibitem{Spivak} M. Spivak, \textit{A Comprehensive Introduction to Differential 
Geometry, Vol V}, 2nd edition, Publish or Perish Inc. (1979).

\bibitem{shankar} S. Venkataramani, \textit{Lower bounds for the energy in a crumpled elastic sheet
-- a minimal ridge},  Nonlinearity  {\bf 17}  (2004),  no. 1, 301--312.
\end{thebibliography}
\end{document}